\documentclass[preprint,3p,12pt,pdf]{elsarticle}

\usepackage{lineno,hyperref}
\modulolinenumbers[5]

\usepackage{epstopdf}

\usepackage{tikz}
\usetikzlibrary{arrows}

\usepackage{pst-node}
\usepackage{tikz-cd}

\usepackage[shellescape]{gmp}

\usepackage{mathrsfs}
\usepackage{amssymb}
\usepackage{amsfonts}
\usepackage{latexsym}
\usepackage{mathtools}
\usepackage{xcolor}
\usepackage{wrapfig}
\usepackage{floatflt}

\usepackage{mathtools}
\usepackage{extarrows}

\usepackage{graphicx}

\usepackage{subcaption}

\usepackage{endnotes}

\usepackage{hyperref}

\def\lb{\label}

\newcommand{\er}[1]{\textrm{(\ref{#1})}}

\newtheorem{theorem}{\bf Theorem}[section]

\def\a{\alpha}         
\def\b{\beta}          
         
\def\G{\Gamma}

           \def\mJ{{\mathscr J}}
        
\def\k{\kappa}

\def\o{\omega}

\def\ve{\varepsilon}       \def\vp{\varphi}    \def\vk{\varkappa}

    \def\R{{\mathbb R}}   \def\C{{\mathbb C}}    
    \def\N{{\mathbb N}}   

\def\lt{\biggl}                  \def\rt{\biggr}
\def\ol{\overline}               \def\wt{\widetilde}

\let\ge\geqslant                 \let\le\leqslant

\def\iy{\infty}
               
\def\ss{\subset}                 
                
                 \def\ev{\equiv}
        
\def\el2{\ell^{\,2}}             \def\1{1\!\!1}

\def\arg{\mathop{\mathrm{arg}}\nolimits}
\def\const{\mathop{\mathrm{const}}\nolimits}

\def\Im{\mathop{\mathrm{Im}}\nolimits}

\def\Re{\mathop{\mathrm{Re}}\nolimits}
\def\Res{\mathop{\mathrm{Res}}\nolimits}

\def\BBox{\hspace{1mm}\vrule height6pt width5.5pt depth0pt \hspace{6pt}}

\let\ge\geqslant
\let\le\leqslant

\newcommand{\ca}{\begin{cases}}
	\newcommand{\ac}{\end{cases}}
\newcommand{\ma}{\begin{pmatrix}}
	\newcommand{\am}{\end{pmatrix}}
\renewcommand{\[}{\begin{equation}}
	\renewcommand{\]}{\end{equation}}
\def\eq{\begin{equation}}
	\def\qe{\end{equation}}

\def\BBox{\hspace{1mm}\vrule height6pt width5.5pt depth0pt \hspace{6pt}}

\begin{document}

	\begin{frontmatter}

		\title{Complete right tail asymptotic for the density of branching processes with fractional generating functions}

		\date{\today}

		\author
		{Anton A. Kutsenko}
	
\address{University of Hamburg, MIN Faculty, Department of Mathematics, 20146 Hamburg, Germany; email: akucenko@gmail.com}

	\begin{abstract}
		The right tail asymptotic series consisting of attenuating exponential terms are derived for the densities of Galton-Watson processes with fractional probability generating functions. 
		The frequencies in the exponential factors form fractal structures in the complex plane.
		We discuss conditions when the asymptotic series converges everywhere. The obtained right tail asymptotic is compared with the standard integral representation of the density and with the complete left tail asymptotic.
	\end{abstract}


	
\end{frontmatter}


{\section{Introduction}\lb{sec0}}

The analysis of the tail asymptotic of the (integral of) density of the martingale limit of supercritical Galton–Watson processes was initiated in \cite{H}. The first asymptotic terms in the right tail were found in \cite{BB1} for the logarithm and clarified without the logarithm in \cite{FW}. In addition to the mentioned papers, one can look at \cite{BD1}, \cite{BD2}, \cite{DM}, \cite{WDK}, and references therein. The current paper considers a special but fairly broad class of rational probability-generating functions. In this case, it is possible to write the complete right tail asymptotic series explicitly. Moreover, sometimes, these series converge to the density of the martingale limit everywhere. We compare the series with two universal formulas for the density: the complete left tail asymptotic series and the Fourier integral representation. Recently, in \cite{FJ}, the universal right-tail bound for the integral of the density is expressed in terms of a simple decaying exponent. Our right-tail asymptotic series also consists of decaying exponents, but their form is quite complex - often, the distribution of frequencies has a fractal structure in the complex plane. Thus, starting from the second asymptotic term (the first main frequency is almost always real), the right tail asymptotic contains a lot of attenuating oscillations. Sometimes, the presence and absence of oscillations in tail asymptotics may attract special attention; see, e.g., \cite{DIL}, \cite{CG}, \cite{DMZ}.

The methods we use involve analytical techniques based on contour integration and residue theory for functions with non-trivial distribution of poles in the complex plane. Generally, these ideas can be used in the Fourier analysis of multiple functional iterations. Such iterative mappings are widely used in practice. A nice introduction to the mathematical theory of iterative mappings is given in \cite{M}. At the same time, the history of such mappings begins with classical works \cite{S1870}, \cite{K1884}, \cite{P1890}, and \cite{F1}.

We consider a simple Galton-Watson branching process in the supercritical case with the minimum family size $1$ - the so-called Schr\"oder case. The Galton--Watson process is defined by
$$
 X_{t+1}=\sum_{j=1}^{X_t}\xi_{j,t},\ \ \ X_0=1,\ \ \ t\in\N\cup\{0\},
$$
where all $\xi_{j,t}$ are independent and identically-distributed natural number-valued random variables with the probability-generating function
$$
 G(z):=\mathbb{E}z^{\xi}.
$$
We assume that it is a ratio of two polynomials
\[\lb{000}
 G(z)=\frac{P(z)}{Q(z)},\ \ P(z)=p_1z+p_2z^2+...+p_Nz^N,\ \ Q(z)=q_0+q_1z+...+q_Mz^M,\ M\ge1.
\]
The probability of the minimum family size is $0<r<1$, where we denote
\[\lb{000a}
 r:=\frac{p_1}{q_0}.
\]
In our case, $p_0=0$. The case of non-zero extinction probability $p_0/q_0\ne0$ can usually be reduced to the supercritical case with the help of Harris-Sevastyanov transformation. The corresponding information about branching processes is given in, e.g., \cite{H} and \cite{B1}. We assume that $P$ and $Q$ have no common factors. All Taylor coefficients at $z=0$ of $G(z)$ are non-negative because they represent the probabilities of the offspring distribution. The generating function satisfies also the identity $G(1)=1$, and 
\[\lb{001}
 E:=G'(1)<+\iy,
\] 
where $E$ is the expectation of the offspring distribution.
Then one may define the {\it martingale limit} $W=\lim_{t\to+\iy}E^{-t}X_t$, the density of which can be expressed as a Fourier transform of some special function
\[\lb{002}
 p(x)=\frac1{2\pi}\int_{-\iy}^{+\iy}\Pi(\mathbf{i}y)e^{-\mathbf{i}yx}dy,
\ \ \ {\rm where}\ \ \ 
 \Pi(z)=\lim_{t\to+\iy}\underbrace{G\circ...\circ G}_{t}(1+\frac{z}{E^{t}}),
\]
see, e.g., \cite{D1}. By definition, it is seen that $\Pi$ satisfies the so-called Poincar\'e-type functional equation
\[\lb{003}
 G(\Pi(z))=\Pi(Ez),\ \ \ \Pi(0)=1,\ \ \ \Pi'(0)=1.
\] 
Since the Fourier integral in \er{002} is quite complex, any expansion and asymptotic analysis of $p(x)$ is welcome. 

\begin{figure}
	\centering
	\begin{subfigure}[b]{0.99\textwidth}
		\includegraphics[width=\textwidth]{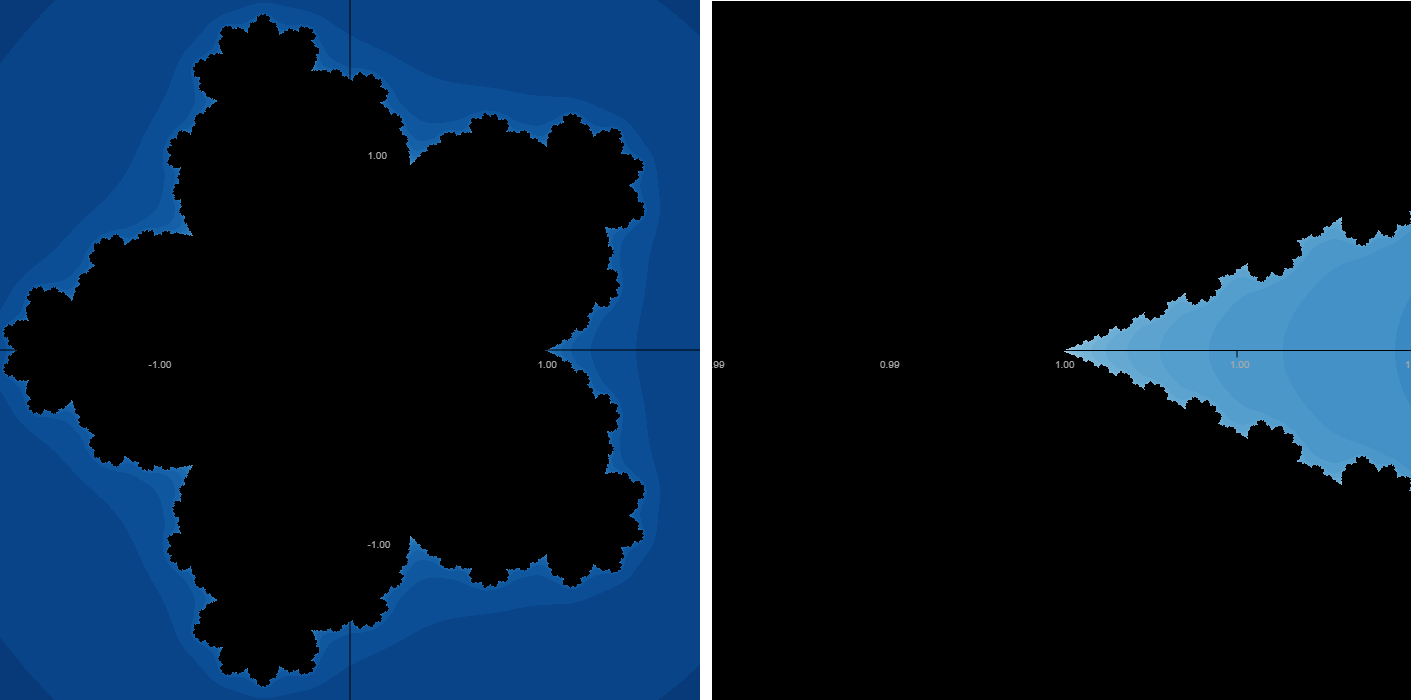}
		\label{fig0a}
	\end{subfigure}
	\hfill
	\begin{subfigure}[b]{0.99\textwidth}
		\includegraphics[width=\textwidth]{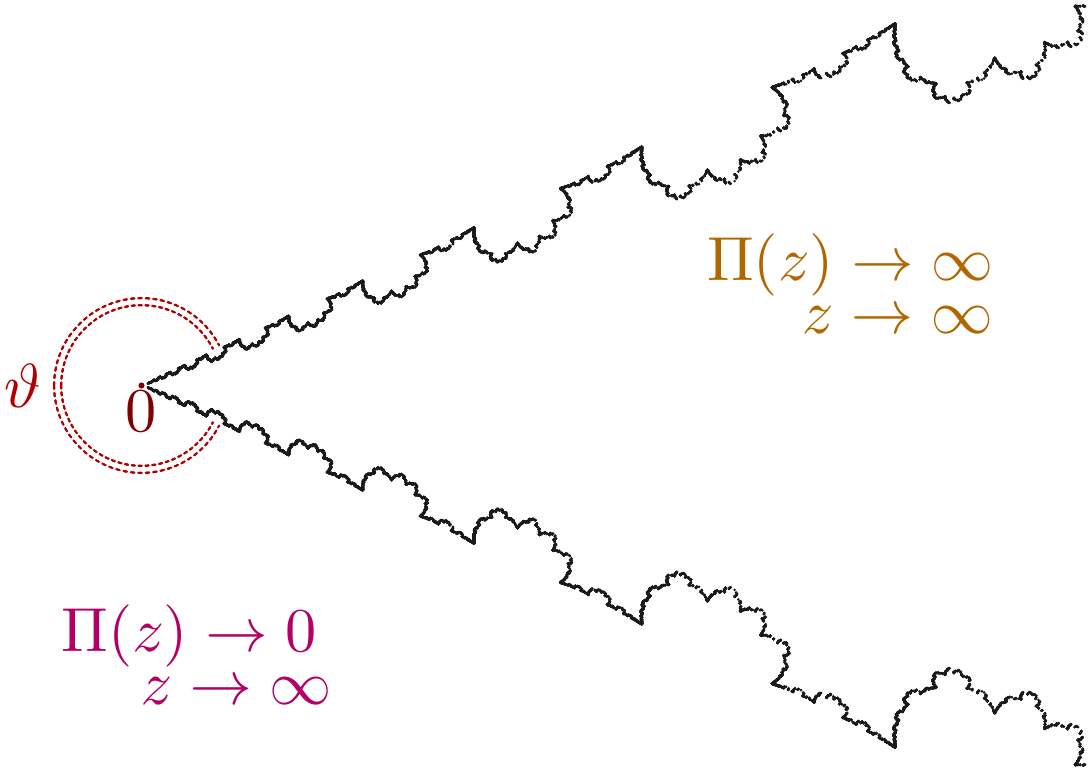}
		\label{fig0b}
	\end{subfigure}
	\caption{The filled Julia set (black area) for $G(z)=0.9z+0.1z^6$ and its zoom at $z=1$. ``Zoomed and shifted'' Julia set divides $\C$ onto two sectors with different behavior of $\Pi(z)$.}\label{fig0}
\end{figure}

To understand what $\Pi(z)$ looks like, we should start with the structure of the Julia set for the generating function $G(z)$. The filled Julia set consists of points $z$ in the complex plane whose orbits are bounded under multiple actions of $G(z)$. Since all the Taylor coefficients of $G(z)$ at $z=0$ are non-negative and $G(1)=1$, we conclude that at least the unit ball $\{z:\ |z|\le1\}$ belongs to the filled Julia set. Moreover, all the points that belong to the open component $\mJ_0$ of the filled Julia set containing the open unit ball satisfy 
\[\lb{004}
 \underbrace{G\circ...\circ G}_{t}(z)\to0\ \  {\rm for}\ \  z\in\mJ_0. 
\] 
This is because $z=0$ is the unique attracting point inside the open unit ball, and $G$ is a contraction mapping here; see details in, e.g., \cite{K24}.
Due to the definition of $\Pi(z)$, see \er{002} and \er{003}, the structure of $\mJ_0$ in the vicinity of $1$ determines the behavior of $\Pi(z)$. Roughly speaking, the ``zoomed'' and transferred to $z=0$ boundary of $\mJ_0$ near $z=1$ divides the complex plane into two parts with different behavior of $\Pi(z)$. This fact is illustrated in Fig. \ref{fig0}. Since $G(z)$ is a polynomial in this example, the right part corresponds to the domain where $\Pi(z)\to\iy$ when $z\to\iy$. If $G(z)$ is a fractional function, the behavior of $\Pi(z)$ on the right part of the complex plane can be different. Moreover, in both cases, the boundary can consist of many connected components. However, the left part always contains the half plane $\{z:\ \Re z<0\}$ since the unit ball belongs to the Julia set. The ``zoomed'' and ``shifted'' unit ball is exactly the left half-plane. One of the main characteristics here is the critical angle $\vartheta$, see Fig. \ref{fig0}, which can be defined as
\[\lb{005}
 \vartheta=\sup\{\vp:\ \{1-re^{\mathbf{i}\vk}\}\ss\mJ_0\ {\rm for}\ |\vk|<\vp\ {\rm and}\ r\to+0\}.
\]
Since $\mJ_0$ contains the unit ball, we have $\vartheta\ge\pi$. For all the examples I tested, $\vartheta>\pi$ strictly. I cannot find in the literature the proof of the statement that $\vartheta>\pi$ for any reasonable $G$, but it seems that it is not difficult to adapt the reasoning from, e.g., \cite{K} involving the second derivative $G''(1)>0$ to obtain the proof. We skip this because it can be cumbersome and, maybe, already done somewhere.

The first two figures in Fig. \ref{fig0} are computed with the help of \href{https://www.marksmath.org/visualization/polynomial_julia_sets/}{this site}\endnote{\it https://www.marksmath.org/visualization/polynomial\_julia\_sets/ \lb{ref1}}. I am not an owner of this resource, but, in my opinion, this is one of the best places on the internet where you can see the visualization of general concepts of holomorphic dynamics. Also, this resource contains utilities for the visualization of Julia sets for rational functions, not only for the polynomials, see \href{https://marksmath.org/visualization/rational_julia_sets/}{this place}\endnote{\it https://www.marksmath.org/visualization/rational\_julia\_sets/ \lb{ref2}}. 

The convergence of orbits \er{004} is exponentially fast, namely, as $r^t$, or faster if one of the orbit points is already $0$, since $G'(0)=r$ and $G(0)=0$, see \er{000} and \er{000a}. This allows us to define an analytical function
\[\lb{006}
 \Phi(z)=\lim_{t\to+\iy}r^{-t}\underbrace{G\circ...\circ G}_{t}(z),\ \ \ z\in\mJ_0.
\] 
It is seen from the definition that $\Phi(z)$ satisfies the so-called Schr\"oder-type functional equation
\[\lb{007}
\Phi(G(z))=r\Phi(z),\ \ \ \Phi(0)=0,\ \ \ \Phi'(0)=1.
\]
The functional equations \er{003} and \er{007} allow us to construct a periodic function very useful in the analysis of $\Pi(z)$. One-periodic Karlin-McGregor function, see, e.g., \cite{KM1}, \cite{KM2}, and \cite{Du1} is defined by
\[\lb{008}
 K(z)=r^{-z}\Phi(\Pi(-E^z)).
\]
Taking into account the periodicity of $K(z)$ and the information about the domains of definition of functions $\Pi(z)$ and $\Phi(z)$ discussed above, it is not difficult to see that $K(z)$ is, at least, defined in the strip $|\Im z|<{\vartheta}/({2\ln E})$. This fact also gives the rate of attenuation of Fourier coefficients 
\[\lb{009}
 K(z)=\sum_{n=-\iy}^{+\iy}\theta_ne^{2\pi\mathbf{i}nz},\ \ |\theta_n|\le C_{\ve}e^{-(\frac{\pi\vartheta}{\ln E}-\ve)|n|}
\]
for any $\ve>0$ and some constants $C_{\ve}>0$. Using \er{008}, we obtain
\[\lb{010}
 \Pi(-z)=\Phi^{-1}(z^{\log_Er}K(\log_Ez)).
\]
Since $K(z)$, as an analytic periodic function, is uniformly bounded in the strip $|\Im z|<(\vartheta-\ve)/(2\ln E)$ for any $\ve>0$, and $\Phi(z)\sim z$ for $z\to0$, from \er{010}, we deduce that
\[\lb{011}
 |\Pi(-z)|\le D_{\ve}|z|^{\log_Er},\ \ \ |\arg z|<\frac{\vartheta-\ve}2
\]
for any $\ve>0$ and some $D_{\ve}>0$. Recall that $\log_Er<0$ because $0<r<1$. Thus, $\Pi(z)$ is decaying to $0$ when $z\to\iy$ is inside the corresponding cone. Everything is ready to write the right tail asymptotic expansion.

\begin{theorem}\lb{T1}
If $\vartheta>\pi$ and $\log_Er<-1$ then
\[\lb{012}
 p(x)\sim -\sum_{\a}\Res(\Pi,\o_{\a})e^{-\o_{\a} x},\ \ \ x\to+\iy,  
\]
where $\{\o_{\a}\}$ are all the poles of $\Pi$, arranged in ascending order of their real parts. Any $\o_{\a}$ can be written in the form
\[\lb{013}
 \o_{\a}=E^{n}\Pi^{-1}(G^{-1}_{i_1}\circ...\circ G^{-1}_{i_m}(z_j)),
\]
where $n>m\ge1$ are some integer numbers, $G_{i_k}^{-1}$ is a branch of multi-valued function $G^{-1}$, and $z_j$ is a zero of $Q(z)$. The residue $\Res(\Pi,\o)$ can be expressed through the values and derivatives of $\Pi$ (and $G$) at points having the form \er{013}. 
\end{theorem} 
{\it Proof.} Suppose $\o$ is a pole of $\Pi$. Then $E^{-n}\o$ is a regular point for some large $n$, since $\Pi(z)$ is analytic for small $z$. By \er{003}, we have 
$$
 \Pi(E^tE^{-n}\o)=\underbrace{G\circ...\circ G}_{t}(E^{-n}\o).
$$
Thus, for some $t<n$ the value $\underbrace{G\circ...\circ G}_{t}(E^{-n}\o)$ should be a pole of $G$, which is a zero of $Q$. Denoting such minimal $t$ as $m$ and the zero as $z_j$, we obtain \er{013}. It is seen also that the residue can be computed in terms of derivatives $\Pi$ and $G$ at the points of the form $\er{013}$. 

Note that if $\deg P>\deg Q$ and $\o$ is a pole of $\Pi$ then $E^s\o$ is a pole of $\Pi$ as well for any $s\in\N$, since $\Pi(Ez)=P(\Pi(z))/Q(\Pi(z))$ and if $\Pi(z)=\iy$ then $\Pi(Ez)=\iy$ as well. This means that $n$ can be greater than $m+1$ in \er{013}. 

The mapping $G^{-1}$ is invertible and contracting in a neighborhood of $z=1$, since $G(1)=1$ and $G'(1)=E>1$. Thus, we can find open neighborhoods of $1$, denote them as $U_1\ss V_1$, such that $G:U_1\to V_1$ and $G^{-1}_0:V_1\to U_1$. Here, $G_0^{-1}$ is the corresponding branch of $G^{-1}$ analytic in $V_1$. The function $\Pi^{-1}$ is also analytic in some open neighborhood of $z=1$, since $\Pi(0)=1$ and $\Pi'(0)=1$. Decreasing, if necessary, the size of $V_1$ (and $U_1$), we may assume that $\Pi: U_0\to V_1$ and $\Pi^{-1}: V_1\to U_0$, where $U_0$ is some open neighborhood of $0$. Thus 
\[\lb{014}
 E\Pi^{-1}(G_0^{-1}(z))=\Pi^{-1}(z),\ \ \  z\in U_1,
\] 
because $\Pi(Ez)=G(\Pi(z))$. Hence, without loss of generality, we may assume that $G^{-1}_{i_1}\circ...\circ G^{-1}_{i_m}(z_j)\not\in U_1$ in \er{013}, since, otherwise $i_1$ should be zero and we can shorten the composition by \er{014}. In turn, this means that $\Pi^{-1}(G^{-1}_{i_1}\circ...\circ G^{-1}_{i_m}(z_j))\not\in U_0$ and, therefore, under the mentioned assumptions $|\o_{\a}|>RE^{n}$ in \er{013} for some constant $R>0$ depending on the size of $U_0$. Since $|\o_{\a}|>RE^{n}$, $n>m$, and the number of branches $G_i^{-1}$ is finite, we deduce that any bounded domain contains only a finite number of the poles. In other words, the set of poles of $\Pi$ has no finite condensation points. Thus, we can choose arbitrary large $t>0$ such that the line $t+\mathbf{i}\R$ does not contain any poles of $\Pi$. Due to \er{002} and the Cauchy's residue theorem, we have
\begin{multline}\lb{015}
 p(x)=\frac1{2\pi\mathbf{i}}\int_{\mathbf{i}\R}\Pi(z)e^{-zx}dz-\frac1{2\pi\mathbf{i}}\int_{t+\mathbf{i}\R}\Pi(z)e^{-zx}dz+\frac1{2\pi\mathbf{i}}\int_{t+\mathbf{i}\R}\Pi(z)e^{-zx}dz
\\
 =-\sum_{\Re\o_{\a}<t}\Res(\Pi,\o_{\a})e^{-\o_{\a} x}+q_t(x),\ \ \ q_t(x):=\frac1{2\pi\mathbf{i}}\int_{t+\mathbf{i}\R}\Pi(z)e^{-zx}dz.
\end{multline}
We have $\vartheta>\pi$ and, hence, both infinite tails of the path $t+\mathbf{i}\R$ must belong to the cone, where \er{011} is fulfilled. Using \er{011}, we obtain
\[\lb{016}
 |q_t(x)|\le e^{-tx}\int_{t+\mathbf{i}\R}D_t|z|^{\log_Er}dy\le \wt{D}_te^{-tx},
\]
for some constants $D_t$ and $\wt{D}_t$, since $\log_Er<-1$ the integral converges.  Identities \er{015} and estimate \er{016} show that asymptotic series \er{012} is true. \BBox

{\bf Remark 1.} While \er{013} is a complete right tail asymptotic series for $p(x)$, in many cases, we cannot replace $\sim$ with $=$. For example, if $Q(z)\ev1$ and $G(z)=P(z)$ is a polynomial of degree $\deg G\ge2$ then $p(x)$ has a super-exponential attenuation of the order $e^{-x^{\b}}$, $x\to+\iy$ with some $\b>1$. The function $\Pi(z)$ is entire and has no poles. It is still attenuating in the left part of $\C$, see Fig. \ref{fig0}. However, due to $\Pi(Ez)=G(\Pi(z))$, it has a super-exponential grows of the order $\exp(z^{\log_E\deg G})$ in the right part of Fig. \ref{fig0}. To obtain the super-exponential attenuation of its Fourier transform $p(x)$ for $x\to+\iy$, one may apply the method of steepest descent. Research on the right tail asymptotic for $p(x)$ in this case is more or less finished in \cite{FW} - at least, for the first asymptotic term.

{\bf Remark 2.} Let us return to the case $Q(z)\ne\const$. If $\deg P\le \deg Q +1$ - we write it as $\deg G \le 1$ - then $\Pi(z)$ has no exponential growth in the right part of the complex plane, see the reasoning set out in the previous remark. Hence, taking arbitrary large $t$ in \er{015}, we can make $q_t(x)$ to be arbitrarily small. Thus, we expect that $\sim$ in \er{012} can be replaced with $=$ in this case.

{\bf Remark 3.} For the Schr\"oder case, research on the left-tail asymptotic was initiated in \cite{D1}, significantly advanced in \cite{BB1}, and finished in \cite{K24}. The corresponding result from \cite{K24} is formulated in terms of the Karlin-McGregor function $K(z)$, see \er{008}.
\begin{theorem}\lb{T2}
If $\vartheta>\pi$ and $\log_Er<-1$ then
\[\lb{017}
p(x)=x^{\a}V_1(x)+x^{\a+\b}V_2(x)+x^{\a+2\b}V_3(x)+...,\ \ \ x>0,
\]
with $\a:=-1-\log_Er>0$, $\b:=-\log_Er>1$, and	
\[\lb{018}
V_m(x)=K_m(\frac{-\ln x}{\ln E}),\ \ \ K_m(z)=\k_m\sum_{n=-\iy}^{+\iy}\frac{\theta_{n}^{\ast m}e^{2\pi\mathbf{i}nz}}{\Gamma(-\frac{2\pi\mathbf{i}n+m\ln r}{\ln E})},\ \ \theta_{n}^{\ast m}=\int_0^1K(x)^me^{-2\pi\mathbf{i}nx}dx,
\]
where $\k_m$ are Taylor coefficients of
\[\lb{019}
 \Phi^{-1}(z)=\k_1z+\k_2z^2+\k_3z^3+....
\]
\end{theorem}
This Theorem is proven in \cite{K24} under the assumption that $G(z)$ is entire, but it is easy to extend the result to rational $G(z)$ and even more big classes of $G(z)$ regular at $z=1$. Using Stirling's approximation, see, e.g., \cite{WW},
$$
 \G(z)=\sqrt{\frac{2\pi}{z}}\lt(\frac{z}{e}\rt)^z\lt(1+O\lt(\frac1z\rt)\rt),\ \ \ |\arg z|<\pi-\ve,
$$
we, for some constant $B>0$, obtain
\[\lb{020}
 \lt|\Gamma\lt(-\frac{2\pi\mathbf{i}n+m\ln r}{\ln E}\rt)\rt|\ge \frac{Be^{-\frac{m\ln r}{\ln E}\ln\sqrt{\frac{4\pi^2n^2+(\ln r)^2m^2}{(\ln E)^2}}-\frac{\pi^2|n|}{\ln E}}}{\sqrt{4\pi^2n^2+(\ln r)^2m^2}}.
\]
The constants $C_{\ve}>0$ in \er{009} can be chosen such that
\[\lb{021}
K(z)^m=\sum_{n=-\iy}^{+\iy}\theta^{\ast m}_ne^{2\pi\mathbf{i}nz},\ \ |\theta^{\ast m}_n|\le (C_{\ve})^me^{-(\frac{\pi\vartheta}{\ln E}-\ve)|n|},\ \ \ \ve>0,\ \ m\ge1.
\]
The standard estimates for Taylor coefficients of analytic functions show $|\k_m|\le K^m$, see \er{019}, for some $K>0$ depending on the convergence radius for $\Phi^{-1}(z)$. Combining these estimates with \er{020} and \er{021}, and using  the conditions $\vartheta>\pi$ and $\ln r<0$ of Theorem \ref{T2}, we obtain
\[\lb{022}
 \lt|\frac{\k_m\theta_{n}^{\ast m}}{\Gamma(-\frac{2\pi\mathbf{i}n+m\ln r}{\ln E})}\rt|\le Ae^{-\a m\ln (m+|n|)-\b|n| }\ \ \ {\rm for\ some}\ \ \ A,\a,\b>0.
\]
Thus, the convergence of the Fourier series in \er{018} is exponentially fast, and the convergence of the series \er{017} is super-exponentially fast.

In the next section, we compare three formulas \er{002}, \er{012}, and \er{017} in a series of examples. We focus on the cases $\deg P\le \deg Q+1$ for which \er{012} should converge to $p(x)$ everywhere $(x>0)$. Also, when $G^{-1}(z)$ admits closed-form expression, \er{012} demonstrates the fastest way to compute $p(x)$. However, if $\deg P> \deg Q+1$ then $\er{012}$ does not coincide with $p(x)$ - the maximum that it gives is an asymptotic expansion for large $x$. At the same time, the other two formulas \er{002} and \er{013} are universal, but the exact behavior of $p(x)$, $x\to+\iy$, is invisible from them.

{\section{Examples}\lb{sec1}}

For most examples, Embarcadero Delphi Rad Studio Community Edition and the library NesLib.Multiprecision is used. This software provides a convenient environment for programming and well-functioning basic functions for high-precision computations. All the algorithms related to the article's subject, including efficient parallelization, are developed by the author (AK). 

{\subsection{General remarks on computational procedures.}\lb{sec1a}} For the numerical results, we need to compute many functions. Let us discuss numerical schemes for some of them. From \er{002}, we have
\[\lb{200}
 \Pi(z)=\lim_{t\to+\iy}\Pi_t(z),\ \ \ \Pi_t(z):=\underbrace{G\circ...\circ G}_{t}(1+\frac{z}{E^{t}}).
\]
Denoting the second Taylor term at $z=1$ as
\[\lb{201}
 G(z)=1+E(z-1)+(z-1)^2G_1(z),
\]
we obtain
\[\lb{202}
 \Pi_{t+1}(z)=\Pi_t\lt(z+\frac{z^2}{E^t}G_1(1+\frac{z}{E^{t+1}})\rt).
\]
Numerical schemes based on \er{202} demonstrate exponentially fast convergence since $E>1$. Differentiating the composition in \er{200}, it is not difficult to obtain a convenient identity for the numerical computation of the derivative $\Pi'(z)$.

From \er{006}, we have
\[\lb{203}
 \Phi(z)=\lim_{t\to+\iy}\Phi_t(z),\ \ \ \Phi_t(z)=r^{-t}\underbrace{G\circ...\circ G}_{t}(z).
\]
Denoting the second Taylor term at $z=0$ as
\[\lb{204}
G(z)=rz+z^2G_0(z),
\]
we obtain
\[\lb{205}
\Phi_{t+1}(z)=\Phi_t(z)+r^{t-1}\Phi_t(z)^2G_0(r^t\Phi_t(z)).
\]
Numerical schemes based on \er{205} demonstrate exponentially fast convergence since $0<r<1$.

From \er{007}, we have
\[\lb{206}
 \Phi^{-1}(rz)=G(\Phi^{-1}(z)).
\]
Substituting Taylor expansions
\[\lb{207}
 G(z)=rz+r_2z^2+r_3z^3+...,\ \ \ \Phi^{-1}(z)=\k_1z+\k_2z^2+\k_3z^3+...
\]
into \er{206} and using initial condition $\k_1=1$, one can express $\k_j$ through $r_j$, step by step. However, $\k_j$ grows exponentially. It makes sense to compute $a^j\k_j$ instead of $\k_j$, where $a$ is small enough. For this, one may set $\k_1=a$ instead of $\k_1=1$.

When $n\to+\iy$, the Fourier coefficients $\theta_{n}^{\ast m}$ should be divided by the small values of the Gamma function, see \er{018} and \er{022}. A good strategy is to compute the normalized Fourier coefficients
\[\lb{208}
 e^{2\pi ny}\theta_{n}^{\ast m}=\int_0^1K(x-\mathbf{i}y)^me^{-2\pi\mathbf{i}nx}dx,\ \ \ n\ge1,\ \ y\ge0,
\]
and use the symmetry $\theta_{-n}^{\ast m}=\ol{\theta_{n}^{\ast m}}$. The parameter $y$ in \er{208} should be chosen close enough to the edge of the strip where $K$ is defined, see above \er{009}. Taking appropriate $a$ and $y$, it is not difficult to achieve accurate and stable computations of the ratio in \er{018}, namely
\[\lb{209}
 \frac{\k_m\theta_{n}^{\ast m}}{\Gamma(-\frac{2\pi\mathbf{i}n+m\ln r}{\ln E})}=\exp\lt(\ln(a^m\k_m)+\ln(e^{2\pi ny}\theta_{n}^{\ast m})-\ln\Gamma(-\frac{2\pi\mathbf{i}n+m\ln r}{\ln E})-m\ln a-2\pi ny\rt). 
\] 
If necessary, one may also try to compute $(bK(z))^m$ with some small $b$ instead of $K(z)^m$, but we didn't have such a need. In addition, note that the computation of $\ln\G$ is much more stable than $\G$.

From \er{003}, we obtain
\[\lb{210}
 \Pi^{-1}(G(z))=E\Pi^{-1}(z),
\]
which is, up to initial conditions and constants, similar to \er{007}. Thus, for the computation of $\Pi^{-1}$, one may use fast algorithms similar to those for $\Phi(z)$, see \er{203}-\er{205}. 

{\subsection{Example 1.}\lb{sec1b}} Let us consider the probability generating function
\[\lb{300}
 G(z)=\frac{z+z^2}{3-z}.
\]
The filled Julia set for $G(z)$ is given in Fig. \ref{fig1}. 

\begin{figure}[h]
	\center{\includegraphics[width=0.5\linewidth]{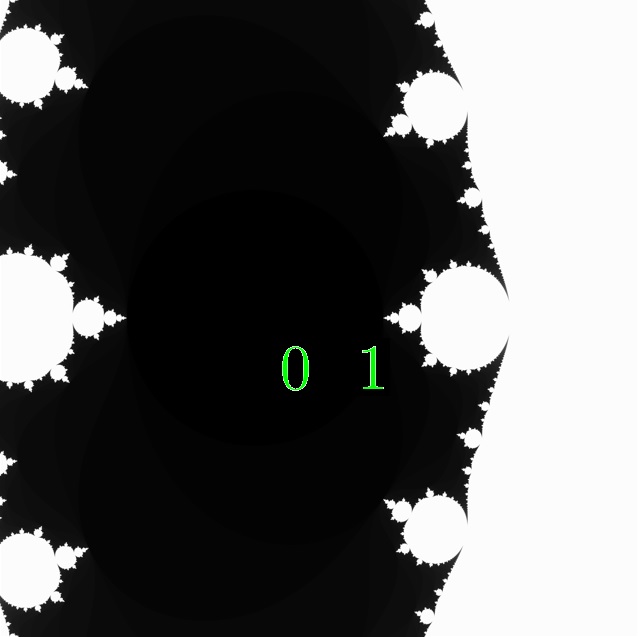}}
	\caption{A fragment of the filled Julia set (black area) for $G$ defined in \er{300}.}\lb{fig1}
\end{figure}

\begin{figure}
	\centering
	\begin{subfigure}[b]{0.45\textwidth}
		\includegraphics[width=\textwidth]{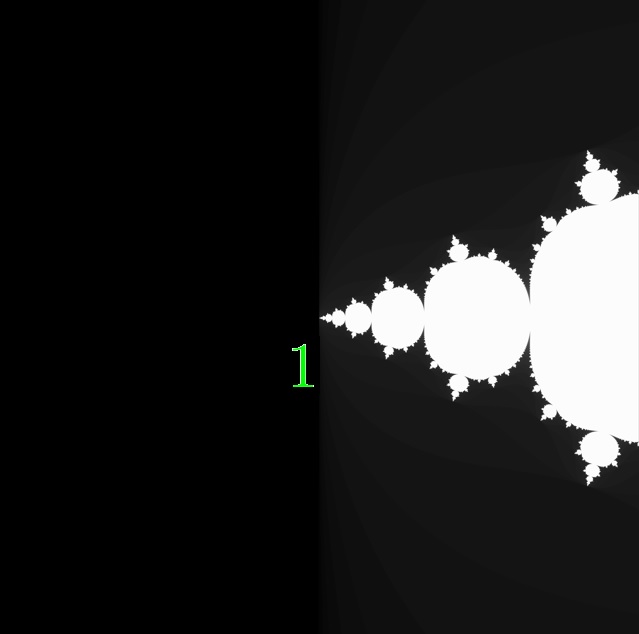}
		\caption{a fragment of the filled Julia set from Fig. \ref{fig1}}
		\label{fig2a}
	\end{subfigure}
	~ 
	\begin{subfigure}[b]{0.45\textwidth}
		\includegraphics[width=\textwidth]{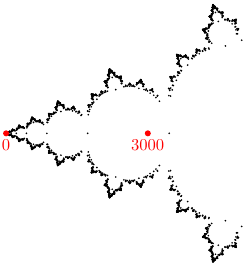}
		\caption{Distribution of poles, see \er{303} and \er{305}}
		\label{fig2b}
	\end{subfigure}
	\caption{(a) The filled Julia set (black area) from Fig. \ref{fig1} near $z=1$ zoomed in about $1000$ times; (b) Distribution of poles of $\Pi(z)$, see \er{002}, in complex plane, for $G(z)$ defined in \er{300}.}\label{fig2}
\end{figure}

There are two branches of the inverse function
\[\lb{301}
 G^{-1}_0(z)=\frac{-1-z+\sqrt{z^2+14z+1}}2,\ \ \ G^{-1}_1(z)=\frac{-1-z-\sqrt{z^2+14z+1}}2.
\]
The probability of the minimal family size and the expectation are
\[\lb{302}
 r=G'(0)=\frac13,\ \ \ E=G'(1)=2.
\]
There is one pole $z=3$ for $G$. We can enumerate all the poles of $\Pi$, see \er{013}, with the help of the following scheme. The first set of ``primary'' poles consists of
\[\lb{303}
 \o_{0,0}=2\Pi^{-1}(3),\ \ \ \o_{j,0}=2^{n+1}\Pi^{-1}\circ G_{j_n}^{-1}\circ...\circ G_{j_1}^{-1}(3),\ \ \ j\ge1,
\]
where
\[\lb{304}
 j=\{j_n...j_1\}_2
\]
is the dyadic representation of $j$. Note that the last element is always $j_n=1$. Hence, we do not count poles twice because $2\Pi^{-1}(G_0^{-1}(z))=\Pi^{-1}(z)$. Due to $\Pi(2z)=G(\Pi(z))$ - if $\Pi(z)=\iy$ then $\Pi(2z)=\iy$ as well, other poles have the form
\[\lb{305}
 \o_{j,k}=2^k\o_{j,0},\ \ \ k\ge1,
\]
again, in accordance with \er{013}. The distribution of $\o_{j,k}$ is plotted in Fig. \ref{fig2}. It resembles the structure of the zoomed Julia set near $z=1$. Such a phenomenon is similar to that of the distribution of zeros of $\Pi(z)$, see, e.g., \cite{K2020}. 

Using $\Pi(2z)=G(\Pi(z))$ with $G$ defined in \er{300}, one can find the expression for the residue
\[\lb{306}
 \Res(\Pi,\o_{j,0})=\frac{-24}{\Pi'(\frac{\o_{j,0}}2)},\ \ \ \Res(\Pi,\o_{j,k})=(-2)^k\Res(\Pi,\o_{j,0}).
\]
Everything is ready to compute $p_a(x)$ by RHS of \er{012}, where we use about $5000$ ``primary'' poles $\o_{j,0}$, $j\le 5000$, and $\o_{j,k}$ generated by the ``primary'' poles, all satisfy $|\o_{j,k}|\le 5000$. For the computation of exact $p(x)$ by \er{002}, we use the interval of integration $y\in[-2\cdot10^5,2\cdot10^5]$ divided by $2\cdot10^7$ points in the trapezoidal rule. For the computation of $p_b(x)$ given by RHS in \er{017}, we use about $30$ terms $V_j(x)$, $j=1,...,30$ for each of which we use about $50$ (or $100$ with the symmetry) Fourier coefficients. For the computation of $\Pi$ and $\Phi$, we use about $150$ iterations. The comparison of $p(x)$, $p_a(x)$, and $p_b(x)$ is given in Fig. \ref{fig3}. They are almost coincide. In this case, the computation of $p_a(x)$ is fastest. However, if $\deg P>\deg Q+1$ in \er{000}, we cannot expect that $p_a(x)$ will coincide with $p(x)=p_b(x)$. In \er{300}, $\deg P=\deg Q+1$ and everything is OK, see {\bf Remark 2} above.

\begin{figure}[h]
	\center{\includegraphics[width=0.8\linewidth]{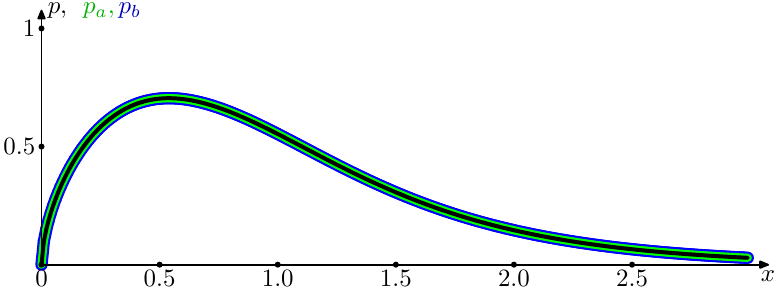}}
	\caption{Comparison of three formulas \er{002}, \er{012}, and \er{017} for the computation of the density $p$ in case of the probability generating function \er{300}.}\lb{fig3}
\end{figure}

Let us discuss the computation of the coefficients $\k_j$, see \er{019}. From \er{206}, we have
\[\lb{307}
 (3-\Phi^{-1}(z))\Phi^{-1}(\frac{z}{3})=\Phi^{-1}(z)+\Phi^{-1}(z)^2,
\]
which leads to
\[\lb{308}
 \k_N=-(3^{N-1}-1)^{-1}\sum_{j=1}^{N-1}\k_j\k_{N-j}(1+3^{-j}),\ \ \ \k_1=1.
\]
Setting $\k_1=a$ with $a=1/5.3$ instead of $\k_1=1$ we compute $a^j\k_j$, which is more accurate than the direct computation of $\k_j$ because of their exponential growth, see the discussion after \er{207}. We use $y=2.9$ in \er{208} to compute the normalized Fourier coefficients. The integrals in \er{208} are computed with $10^6$ uniformly distributed nodes in the trapezoidal rule.

{\subsection{Example 2.}\lb{sec1c}}

Let us consider the probability-generating function
\[\lb{400}
 G(z)=\frac{4(z+z^2)}{9-z^2}.
\]
The filled Julia set for $G(z)$ is given in Fig. \ref{fig4}.

\begin{figure}[h]
	\center{\includegraphics[width=0.5\linewidth]{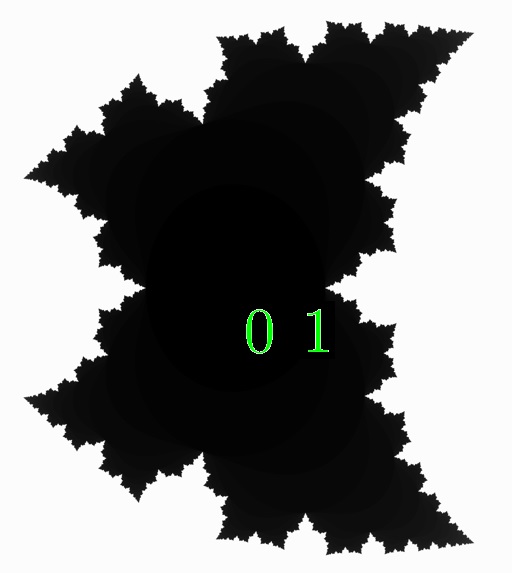}}
	\caption{The filled Julia set for $G$ defined in \er{400}.}\lb{fig4}
\end{figure}

\begin{figure}
	\centering
	\begin{subfigure}[b]{0.45\textwidth}
		\includegraphics[width=\textwidth]{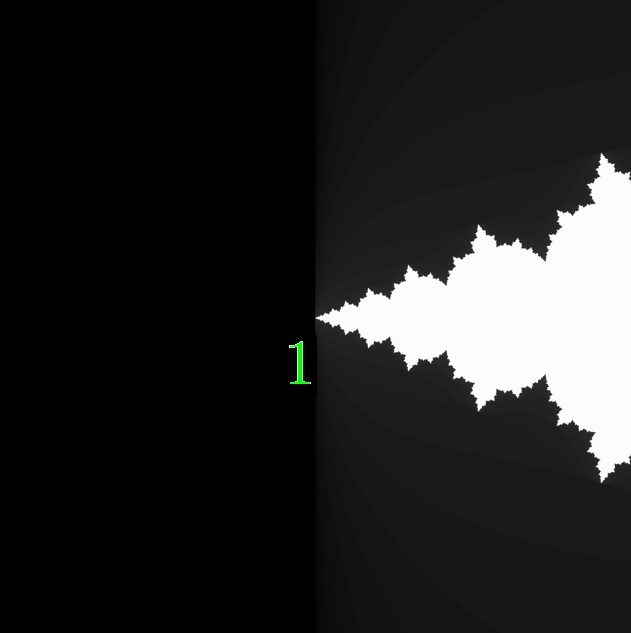}
		\caption{a fragment of the filled Julia set from Fig. \ref{fig4}}
		\label{fig5a}
	\end{subfigure}
	~ 
	\begin{subfigure}[b]{0.45\textwidth}
		\includegraphics[width=\textwidth]{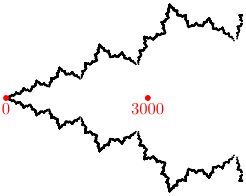}
		\caption{Distribution of poles, see \er{403} and \er{404}}
		\label{fig5b}
	\end{subfigure}
	\caption{(a) The filled Julia set (black area) from Fig. \ref{fig4} near $z=1$ zoomed in about $1000$ times; (b) Distribution of poles of $\Pi(z)$, see \er{002}, in complex plane, for $G(z)$ defined in \er{400}.}\label{fig5}
\end{figure}

There are two branches of the inverse function
\[\lb{401}
 G^{-1}_0(z)=\frac{-2+\sqrt{9z(z+4)+4}}{4+z},\ \ \ G^{-1}_1(z)=\frac{-2-\sqrt{9z(z+4)+4}}{4+z}.
\]
The probability of the minimal family size and the expectation are
\[\lb{402}
 r=G'(0)=\frac49,\ \ \ E=G'(1)=\frac74.
\]
There are two poles $z=3$ and $z=-3$ for $G$. We can enumerate all the poles of $\Pi$, see \er{013}, with the help of the following scheme. The first set of poles consists of
\[\lb{403}
 \o_{0,0}=\frac74\Pi^{-1}(3),\ \ \ \o_{j,0}=\lt(\frac74\rt)^{n+1}\Pi^{-1}\circ G_{j_n}^{-1}\circ...\circ G_{j_1}^{-1}(3),\ \ \ j\ge1,
\]
and the second set
\[\lb{404}
 \o_{0,1}=\frac74\Pi^{-1}(-3),\ \ \ \o_{j,1}=\lt(\frac74\rt)^{n+1}\Pi^{-1}\circ G_{j_n}^{-1}\circ...\circ G_{j_1}^{-1}(-3),\ \ \ j\ge1,
\]
where
\[\lb{405}
 j=\{j_n...j_1\}_2
\]
is the dyadic representation of $j$. 
In contrast to the previous example, if $\o$ is a pole of $\Pi$, then $7\o/4$ is not a pole because $G(\iy)\ne \iy$. The distribution of $\o_{j,k}$ is plotted in Fig. \ref{fig5}. It resembles the structure of the zoomed Julia set near $z=1$.

Using $\Pi(7z/4)=G(\Pi(z))$ with $G$ defined in \er{400}, one can find the expression for the residue
\[\lb{406}
\Res(\Pi,\o_{j,0})=\frac{-14}{\Pi'(\frac{4\o_{j,0}}7)},\ \ \ \Res(\Pi,\o_{j,1})=\frac{7}{\Pi'(\frac{4\o_{j,1}}7)}.
\]
While the distribution of the poles has some certain structure, the distribution of residues is mostly chaotic. We plot some of them in Fig. \ref{fig5r}.

\begin{figure}[h]
	\center{\includegraphics[width=0.7\linewidth]{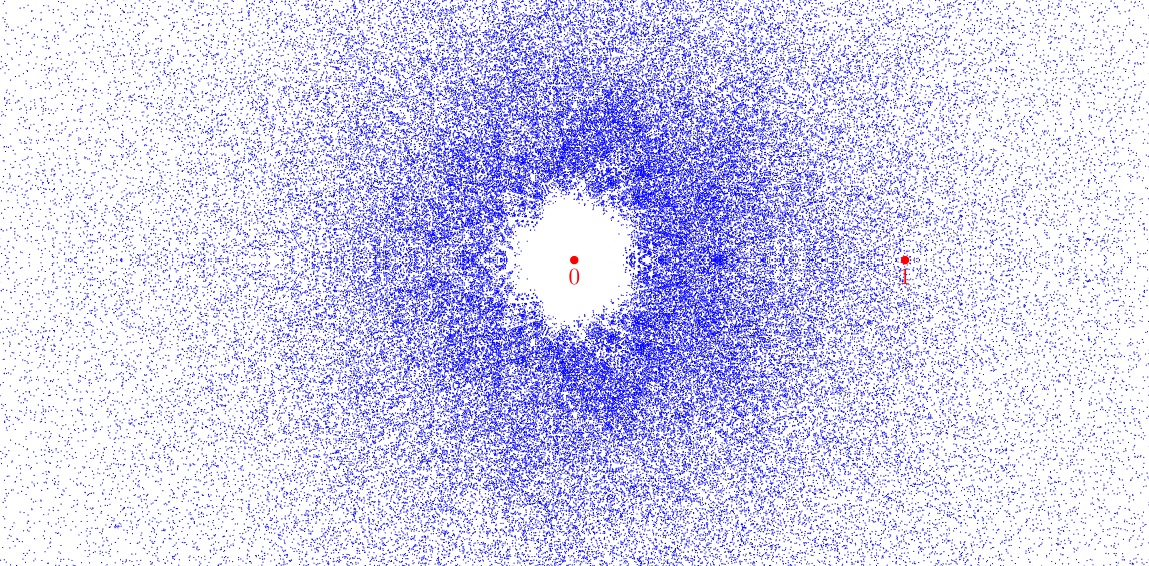}}
	\caption{First $14\cdot10^4$ residues \er{406} are plotted in the complex plane. 
		}\lb{fig5r}
\end{figure}

Everything is ready to compute $p_a(x)$ by RHS of \er{012}, where we use about $2500$ poles $\o_{j,0}$, $j\le 2500$, and the same amount of $\o_{j,1}$: among them we select only those that satisfy $|\o_{j,k}|\le 5000$.  For the computation of exact $p(x)$ by \er{002}, we use the interval of integration $y\in[-2\cdot10^5,2\cdot10^5]$ divided by $2\cdot10^7$ points in the trapezoidal rule. For the computation of $p_b(x)$ given by RHS in \er{017}, we use about $40$ terms $V_j(x)$, $j=1,...,40$ for each of which we use about $60$ (or $120$ with the symmetry) Fourier coefficients. In comparison to the previous example, we should increase the number of harmonics to achieve the same accuracy about $10^{-4}$ on the interval $x\in(0,5)$. If we take $x\in(0,3)$, as it is plotted in Fig. \ref{fig6}, then a much smaller number of harmonics is acceptable. For the computation of $\Pi$ and $\Phi$, we use about $150$ iterations. The comparison of $p(x)$, $p_a(x)$, and $p_b(x)$ is given in Fig. \ref{fig6}. They are almost identical. In this case, the computation of $p_a(x)$ is fastest. However, if $\deg P>\deg Q+1$ in \er{000}, we cannot expect that $p_a(x)$ will coincide with $p(x)=p_b(x)$. In \er{400}, $\deg P=\deg Q$ and everything is OK, see {\bf Remark 2} above.

\begin{figure}[h]
	\center{\includegraphics[width=0.8\linewidth]{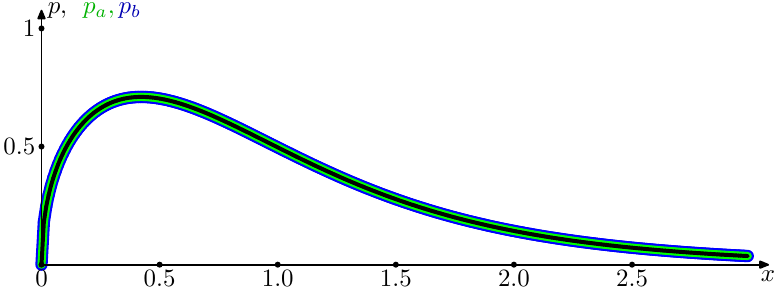}}
	\caption{Comparison of three formulas \er{002}, \er{012}, and \er{017} for the computation of the density $p$ in case of the probability generating function \er{400}.}\lb{fig6}
\end{figure}

Let us discuss the computation of the coefficients $\k_j$, see \er{019}. From \er{206}, we have
\[\lb{407}
(9-\Phi^{-1}(z)^2)\Phi^{-1}(\frac{4z}{9})=4(\Phi^{-1}(z)+\Phi^{-1}(z)^2),
\]
which leads to
\[\lb{408}
\k_N=\lt(4\lt(\frac{4}{9}\rt)^{N-1}-4\rt)^{-1}\lt(4\sum_{j=1}^{N-1}\k_j\k_{N-j}+\sum_{j=1}^{N-2}\lt(\frac{4}{9}\rt)^j\k_j\sum_{k=1}^{N-j-1}\k_k\k_{N-j-k}\rt),\ \ \ \k_1=1,
\]
where the second double sum is zero for $N=2$. Setting $\k_1=a$ with $a=1/4.1$ instead of $\k_1=1$, we compute $a^j\k_j$, which is more accurate than the direct computation of $\k_j$ because of their exponential growth, see the discussion below \er{207}. We use $y=4$ in \er{208} to compute the normalized Fourier coefficients. The integrals in \er{208} are computed with $10^6$ uniformly distributed nodes in the trapezoidal rule.

{\subsection{Example 3.}\lb{sec1d}}

Let us consider the probability-generating function
\[\lb{500}
G(z)=\frac{2z}{(3-z)(2-z)}.
\]
The (filled) Julia set for $G$ is almost the whole complex plane except for some sparse real sets. There are two branches of the inverse function
\[\lb{501}
G^{-1}_0(z)=\frac{5z+2-\sqrt{z^2+20z+4}}{2z},\ \ \ G^{-1}_1(z)=\frac{5z+2+\sqrt{z^2+20z+4}}{2z}.
\]
The probability of the minimal family size and the expectation are
\[\lb{502}
r=G'(0)=\frac13,\ \ \ E=G'(1)=\frac52.
\]
There are two poles $z=2$ and $z=3$ for $G$. We can enumerate all the poles of $\Pi$, see \er{013}, with the help of the following scheme. The first set of poles consists of
\[\lb{503}
\o_{0,0}=\frac52\Pi^{-1}(2),\ \ \ \o_{j,0}=\lt(\frac52\rt)^{n+1}\Pi^{-1}\circ G_{j_n}^{-1}\circ...\circ G_{j_1}^{-1}(2),\ \ \ j\ge1,
\]
and the second set
\[\lb{504}
\o_{0,1}=\frac52\Pi^{-1}(3),\ \ \ \o_{j,1}=\lt(\frac52\rt)^{n+1}\Pi^{-1}\circ G_{j_n}^{-1}\circ...\circ G_{j_1}^{-1}(3),\ \ \ j\ge1,
\]
where
\[\lb{505}
j=\{j_n...j_1\}_2
\]
is the dyadic representation of $j$. Using $\Pi(5z/2)=G(\Pi(z))$ with $G$ defined in \er{500}, one can find the expression for the residue
\[\lb{506}
\Res(\Pi,\o_{j,0})=\frac{-10}{\Pi'(\frac{2\o_{j,0}}5)},\ \ \ \Res(\Pi,\o_{j,1})=\frac{15}{\Pi'(\frac{2\o_{j,1}}5)}.
\]
The poles \er{503}, \er{504} and the residues \er{506} are all real. Hence, in contrast to the previous two examples, there are no oscillations in representation \er{012}. Everything is ready to compute $p_a(x)$ by RHS of \er{012}, where we use about $2500$ poles $\o_{j,0}$, $j\le 2500$, and the same amount of $\o_{j,1}$.  For the computation of exact $p(x)$ by \er{002}, we use the interval of integration $y\in[-2\cdot10^5,2\cdot10^5]$ divided by $2\cdot10^7$ points in the trapezoidal rule. For the computation of $p_b(x)$ given by RHS in \er{017}, we use about $30$ terms $V_j(x)$, $j=1,...,30$ for each of which we use about $50$ (or $100$ with the symmetry) Fourier coefficients. For the computation of $\Pi$ and $\Phi$, we use about $150$ iterations. The comparison of $p(x)$, $p_a(x)$, and $p_b(x)$ is given in Fig. \ref{fig6}. They are almost identical. In this case, the computation of $p_a(x)$ is fastest. However, if $\deg P>\deg Q+1$ in \er{000}, we cannot expect that $p_a(x)$ will coincide with $p(x)=p_b(x)$. In \er{500}, $\deg P=\deg Q-1$ and everything is OK, see {\bf Remark 2} above.

\begin{figure}[h]
	\center{\includegraphics[width=0.8\linewidth]{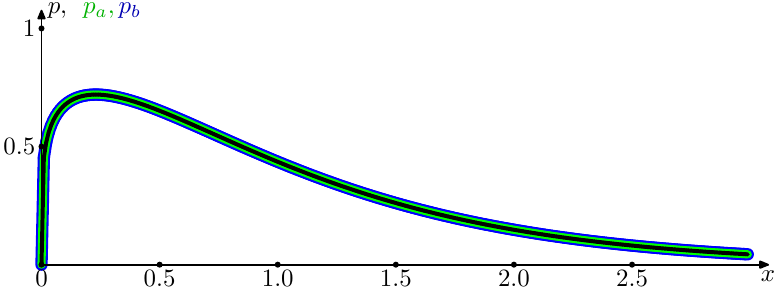}}
	\caption{Comparison of three formulas \er{002}, \er{012}, and \er{017} for the computation of the density $p$ in case of the probability generating function \er{500}.}\lb{fig7}
\end{figure}

Let us discuss the computation of the coefficients $\k_j$, see \er{019}. From \er{206}, we have
\[\lb{507}
(3-\Phi^{-1}(z))(2-\Phi^{-1}(z))\Phi^{-1}(\frac{z}{3})=2\Phi^{-1}(z),
\]
which leads to
\[\lb{508}
\k_N=\lt(2-\frac2{3^{N-1}}\rt)^{-1}\lt(-5\sum_{j=1}^{N-1}\frac{\k_j\k_{N-j}}{3^j}+\sum_{j=1}^{N-2}\frac{\k_j}{3^j}\sum_{k=1}^{N-j-1}\k_k\k_{N-j-k}\rt),\ \ \ \k_1=1,
\]
where the second double sum is zero if $N=2$. Setting $\k_1=a$ with $a=1/2$ instead of $\k_1=1$ we compute $a^j\k_j$, which is more accurate than the direct computation of $\k_j$ because of their exponential growth, see the discussion after \er{207}. We use $y=3$ in \er{208} to compute the normalized Fourier coefficients. The integrals in \er{208} are computed with $10^6$ uniformly distributed nodes in the trapezoidal rule.

\section*{Acknowledgements} 
This paper is a contribution to the project S1 of the Collaborative Research Centre TRR 181 "Energy Transfer in Atmosphere and Ocean" funded by the Deutsche Forschungsgemeinschaft (DFG, German Research Foundation) - Projektnummer 274762653. 

%

\bibliographystyle{abbrv}

\theendnotes

\end{document}